\documentclass[11pt]{article}
\usepackage{amsmath,amsthm}
\usepackage{bm,amssymb}
\usepackage{setspace}
\usepackage{pdfsync}

\newcommand{\1}{\mathbf{1}}

\renewcommand{\P}{\mathrm{P}}
\newcommand{\E}{\mathrm{E}}
\newcommand{\R}{\mathbf{R}}
\renewcommand{\d}{\delta}

\renewcommand{\d}{\mathrm{d}}
\newcommand{\e}{\mathrm{e}}
\renewcommand{\Re}{\mathrm{Re}}

\newtheorem{stat}{Statement}[section]
\newtheorem{proposition}[stat]{Proposition}

\newtheorem{theorem}[stat]{Theorem}
\newtheorem{lemma}[stat]{Lemma}
\theoremstyle{definition} 

\newtheorem{remark}[stat]{Remark}
\newtheorem*{OP}{Open Problems}

\numberwithin{equation}{section}

\begin{document}\onehalfspacing

\title{\bf On the existence and position of the farthest
	peaks of a family of stochastic heat and wave equations%
	\thanks{Research supported in part by the
	Swiss National Science Foundation Fellowship PBELP2-122879 (D.C.)
	and the NSF grant DMS-0706728 (D.K.).}
}
\author{Daniel Conus\\University of Utah
\and 
Davar Khoshnevisan\\University of Utah
}

\date{September 28, 2010}
\maketitle

\begin{abstract}
	 We study the stochastic heat equation
   $\partial_t u = \mathcal{L}u+\sigma(u)\dot W$
   in $(1+1)$ dimensions, where $\dot W$ is space-time white noise,
   $\sigma:\R\to\R$ is Lipschitz continuous, and $\mathcal{L}$ is the generator
   of a symmetric L\'evy process that 
   has finite exponential moments, and $u_0$ has exponential decay at $\pm\infty$.
   We prove that under natural conditions on $\sigma$: 
   (i) The $\nu$th absolute
   moment of the solution to our stochastic heat
   equation grows exponentially with time; and (ii) The distances to the origin
   of the farthest high peaks of those moments grow exactly linearly with time. Very little
   else seems to be known about the location of the high peaks of
   the solution to the stochastic heat equation under the present setting. 
   [See, however, \cite{GM2,GKM} for the analysis of the location of the peaks in a different
   model.]
   
   Finally, we show that these results extend to the stochastic wave equation driven by Laplacian. \\

   \smallskip
	
\end{abstract}
\newpage
\section{Introduction}

We study the nonlinear stochastic heat equation
\begin{equation}\label{heat}
	\frac{\partial}{\partial t} u_t(x) = (\mathcal{L}u_t)(x) +
	\sigma(u_t(x))\frac{\partial^2}{\partial t\partial x} W(t\,,x)\qquad \text{for }
	t>0,\, x\in\R,
\end{equation}
where: (i) $\mathcal{L}$ is the generator of a real-valued symmetric L\'evy process 
$\{X_t\}_{t \geq 0}$ with L\'evy exponent $\Psi$;\footnote{Recall that
$\Psi$ is defined by
$\E\exp(i\xi X_1)=\exp(-\Psi(\xi))$ $[\xi\in\R]$. Because of the
symmetry of $\{X_t\}_{t\ge 0}$,
$\Psi(\xi)=\Psi(-\xi)\ge 0$ for all $\xi\in\R$.}
(ii) $\sigma:\R\to\R$ is Lipschitz continuous with Lipschitz constant 
${\rm Lip}_\sigma$; (iii) $W$ is two-parameter Brownian sheet, indexed by
$(t\,,x)\in\R_+\times\R$; and (iv) the initial function $u_0:\R\to\R_+$ is in $L^\infty(\R)$.
Equation \eqref{heat} arises  for several reasons that include its connections to the
stochastic Burger's equation (see Gy\"ongy and Nualart \cite{GN})
and the parabolic Anderson model (see Carmona and Molchanov \cite{CM94}).

According to the theory of Dalang \cite{Dalang}, \eqref{heat} has a unique
solution when
\begin{equation}\label{Dalang}
	\Upsilon(\beta):= 
	\frac{1}{2\pi}
	\int_{-\infty}^\infty \frac{\d\xi}{\beta+2\Psi(\xi)}<\infty
	\quad\text{for some, hence all, $\beta>0$}.
\end{equation}
Moreover, under various conditions on $\sigma$,
\eqref{Dalang} is necessary for the existence of
a solution \cite{Dalang,PZ}. 

Foondun and Khoshnevisan \cite{FKinter}
have shown that:
\begin{equation}\label{eq:UB}
	\overline\gamma(\nu):=
	\limsup_{t\to\infty} \frac1t \sup_{x\in\R}\ln\E\left(|u_t(x)|^{\nu}\right)<\infty
	\qquad\text{for every $\nu\ge 2$};
\end{equation}
and that
\begin{equation}
	\limsup_{t\to\infty}\frac1t\inf_{x\in\R}\ln\E\left( |u_t(x)|^{\nu}\right)>0
	\qquad\text{for every $\nu \geq 2$,}
\end{equation}
provided that: (a) $\inf_x|\sigma(x)/x|>0$;
and (b) $\inf_x u_0(x)>0$.\footnote{In fact, these results do not require
that $\{X_t\}_{t\ge 0}$ is a symmetric process provided that we replace $\Psi$
with $\Re\Psi$ in \eqref{Dalang}.} 
Together these results show that if $u_0$ is bounded 
away from 0 and $\sigma$ is sublinear,
then the solution to \eqref{heat} is ``weakly intermittent''
[that is, highly peaked for large $t$].  Rather than describe why this
is a noteworthy property, we refer the interested reader to
the extensive bibliography of \cite{FKinter}, which contains several pointers to the literature in mathematical physics  that motivate [weak] intermittency.

The case that $u_0$ has compact support arises equally naturally in mathematical 
physics, but little is known rigorously about when, why, or if the solution to
\eqref{heat} is weakly intermittent when $u_0$ has compact support. 
In fact, we know of only one article 
\cite{FKcpt}, which considers the special case
$\mathcal{L}=\partial^2/\partial x^2$, $\sigma(0)=0$, and $u_0$ 
smooth and compactly supported. 
In that article it is shown that $\overline\gamma(2)\in(0\,,\infty)$, but the
arguments of \cite{FKcpt} rely critically on several special properties of 
the Laplacian. A closely-related case ($u_0:=\delta_0$)
appears in Bertini and Cancrini \cite{BC}.

Presently, we show that weak intermittency follows in some cases 
from a ``stochastic weighted Young inequality''
(Proposition \ref{pr:young:1}).
Such an inequality is likely to have other applications as well.
And more significantly, we describe quite precisely the location of the 
high peaks that are farthest away from the origin.

From now on, let us assume further that
\begin{equation}\label{cond:sigma}
	\sigma(0)=0
	\quad\text{and}\quad
	{\rm L}_\sigma:=\inf_{x\in\R}|\sigma(x)/x|>0.
\end{equation}
And we define two \emph{growth indices}:
\begin{equation} \label{def:lambda_up}
	\overline{\lambda}(\nu) :=\inf\left\{\alpha>0:\
	\limsup_{t\to\infty} \frac 1t \sup_{|x|\ge\alpha t}\ln
	\E\left(|u_t(x)|^{\nu}\right)<0\right\};
\end{equation}
where $\inf\varnothing:=\infty$; and
\begin{equation} \label{def:lambda_down}
	\underline{\lambda}(\nu) :=\sup\left\{\alpha>0:\
	\limsup_{t\to\infty} \frac 1t \sup_{|x|\ge\alpha t}\ln
	\E\left(|u_t(x)|^{\nu}\right)>0\right\};
\end{equation}
where $\sup\varnothing:=0$. 

One can check directly that
$0\le \underline\lambda(\nu)\le \overline\lambda(\nu)\le\infty.$
Our goal is to identify several instances when
$0< \underline\lambda(\nu)\le \overline\lambda(\nu)<\infty.$
In those instances, it follows that: (i) The solution to 
\eqref{heat} has very high peaks as $t\to\infty$
[``weak intermittency'']; and 
(ii) The distances between the origin and
the  farthest high peaks grow exactly linearly in $t$.
This seems to be the first concrete piece of information on
the location of the high peaks of the solution to \eqref{heat}
when $u_0$ has compact support.

Let $\mathcal{D}_{\text{\it exp}}$ denote the collection of all bounded 
lower semicontinuous
functions $h:\R\to\R_+$ for which there exists $\rho>0$ such that $h(x)=O(\e^{-\rho|x|})$
as $|x|\to\infty$.

\begin{theorem}\label{th:main}
	If there exists $c>0$ such that $\E[\e^{cX_1}] < \infty$
	and $u_0\in\mathcal{D}_{\text{\it exp}}$ is strictly positive on
	a set of positive measure, then 
	$0<\underline{\lambda}(\nu)\le \overline{\lambda}(\nu)<\infty$
	for all \ $\nu\in[2\,,\infty)$.
\end{theorem} 

\begin{remark}\label{rem:Trunc:Stable}
	Theorem \ref{th:main} applies to many
	L\'evy processes other than Brownian motion. Here we mention a
	simple family of examples. First, let us recall the L\'evy--Khintchine
	formula for $\Psi$: There exists $\sigma\in\R$
	and a symmetric Borel measure $m$ on $\R$ such that
	$m(\{0\})=0$, $\int_{-\infty}^\infty (1\wedge \xi^2)\, m(\d\xi)<\infty$,
	and for all $\xi\in\R$,
	\begin{equation}
		\Psi(\xi) =\sigma^2\xi^2+2\int_{-\infty}^\infty
		\left[ 1-\cos(\xi z)\right]\, m(\d z).
	\end{equation}
	It is well known that 
	for all rapidly-decreasing functions $f:\R\to\R$,
	\begin{equation}
		(\mathcal{L}f)(x)
		=\sigma^2 f''(x) + \int_0^\infty\left[ f(x+z)+f(x-z)-2f(x)
		\right]\, m(\d z).
	\end{equation}
	It is possible to show that the conditions of Theorem \ref{th:main} are
	met, for example, if the support of $m$ is bounded.
	One can frequently verify Dalang's condition \eqref{Dalang} in such examples,
	as well. For instance, let us consider the particular case that $X$ is a
	``truncated symmetric stable'' process. That is the case when $\sigma:=0$ and
	the L\'evy measure satisfies $m(\d z)=|z|^{-(1+\alpha)}\1_{(-1,1)}(z)\, \d z$
	with $1<\alpha<2$. In this case,
	\begin{equation}
		(\mathcal{L}f)(x) = \int_0^1
		\left[ \frac{f(x+z)-f(x-z)-2f(x)}{z^{1+\alpha}}\right]\d z,
	\end{equation}
	and \eqref{Dalang} holds because
	\begin{equation}
		\Psi(\xi) = 2\int_0^1\frac{1-\cos(\xi z)}{z^{1+\alpha}}\, \d z 
		= (2+o(1))\int_0^\infty\frac{1-\cos r}{r^{1+\alpha}}\, \d r
		\cdot  |\xi|^\alpha,
	\end{equation}
	as $|\xi|\to\infty$.
	More interesting examples can be found within the 
	constructions of Rosi\'{n}ski \cite{Rosinski}
	and Houdr\'e and Kawai \cite{HK}. 

\end{remark}

There are concrete instances where one can improve
the results of Theorem \ref{th:main}, thereby
establish quite good estimates for 
$\underline{\lambda}(2)$ and $\overline{\lambda}(2)$.
The following typifies a good example, in which $\mathcal{L}$
is a constant multiple of the Laplacian.

\begin{theorem} \label{th:heat_exact}
	If $\mathcal{L}f = \frac\kappa 2 f''$ 
	and $u_0$ is lower semicontinuous and has a compact support of positive measure, 
	then Theorem \ref{th:main} holds. In addition,
	\begin{equation} \label{heat_exact_cond}
		\frac{{\rm L}_\sigma^2}{2\pi} \le \underline{\lambda}(2) 
		\le \overline{\lambda}(2) \le \frac{{\rm Lip}_\sigma^2}{2}
		\qquad\text{for all $\kappa>0$}.
	\end{equation}
\end{theorem}

In the case of the Parabolic Anderson Model [$\sigma(u) := \lambda u$],
\eqref{heat_exact_cond} tells us that 
$ \lambda^2/2\pi \le \underline{\lambda}(2) \le \overline{\lambda}(2) \le \lambda^2/2$.

We know from Theorem \ref{th:main} that the positions of the farthest peaks
grow linearly with time. Theorem \ref{th:heat_exact} describes an explicit
interval in which the farthest high peaks necessarily fall. Moreover,
this interval does not depend on the value of the diffusion coefficient
$\kappa$. In intuitive terms, these remarks can be summed up as follows:
``Any amount of noise leads to totally intermittent behavior.''
This observation was made, much earlier, in various physical contexts; see, for example,
Zeldovich, Ruzmaikin, and Sokoloff \cite[pp.\ 35--37]{ZRS}.

We mention that the main ideas in the proofs of Theorems \ref{th:main} 
and \ref{th:heat_exact} apply also in other settings. 
For example, in Section \ref{sec:wave} below
we study a hyperbolic SPDE, and prove that 
$\underline{\lambda}(2)=\underline{\lambda}(\nu) = 
\overline{\lambda}(\nu)=\overline{\lambda}(2)$ for $\nu \ge 2$,
under some regularity hypotheses.
This implies the existence of a sharp phase transition between
exponential growth and exponential decay of those hyperbolic SPDEs. Moreover, we will see that
the intermittent behavior of the stochastic wave equation differs 
from \eqref{heat} in two fundamental ways:
(a) The variance of the noise affects the strength of intermittency; and
(b) the rate of growth of $\sigma$ does not.

We conclude the introduction with two questions that have eluded us.

\begin{OP}\begin{enumerate}
	\item Is there a unique phase transition in the exponential growth of
		\eqref{heat}. In other words, we ask:
		\[
			\text{Is }\ \underline\lambda(\nu)=\overline\lambda(\nu)?
		\]
		Although we have no conjectures about this in the present
		setting of parabolic equations, 
		Theorem \ref{thm:wave} below answers this question 
		affirmatively for some hyperbolic SPDEs.
	\item Suppose $u_0\in\mathcal{D}_{\it exp}$
		and $\mathcal{L}=-(-\Delta)^{\alpha/2}$ denote
		the fractional Laplacian
		for some exponent $\alpha\in(1\,,2)$.
		Does $\sup_{x\in\R}\E(|u_t(x)|^2)$ grow exponentially with $t$?
		We mention the following related fact:
		It is possible to adapt the proof of  \cite[Theorem 2.1]{FKcpt}
		to show that if $u_0\in L^2(\R)$,
		then $\int_{-\infty}^\infty \E(|u_t(x)|^2)\,\d x$ grows exponentially
		with $t$. The remaining difficulty is to establish ``localization.''
		The results of the present paper accomplish all this if
		the fractional Laplacian---which is the generator of a symmetric
		stable process---were replaced by the generator of a truncated
		symmetric stable process; see Remark \ref{rem:Trunc:Stable}.
	\end{enumerate}
\end{OP}

Before proceeding to the proofs of Theorems \ref{th:main}  and \ref{th:heat_exact}, we introduce some notation. We write $\|\cdot\|_{\nu}$ the standard norm on $L^{\nu}(\P)$. That is,
$$\|Y\|_{\nu}:=\{\E(|Y|^{\nu})\}^{1/\nu}, \qquad\text{for all $\nu\in[1\,,\infty)$ and $Y\in L^{\nu}(\P)$}.$$ We now recall the following form of
Burkholder's inequality that will be used here and throughout.

\begin{theorem}[The Burkholder--Davis--Gundy 
	inequality \cite{Burkholder,BDG,BG}]
	Let $\{M_t\}_{t \geq 0}$ be a continuous
	martingale. Then, for all $k \geq 1$ and for all $t > 0$
	there exists a constant $z_k$ such that
	\begin{equation}
		\| M_t  \|_k\le 
		z_k  \|  \langle M \rangle_t  \|_{k/2}^{1/2},
	\end{equation}
	where $\langle M \rangle$ denotes the quadratic variation of $M$.
	\label{th:BDG}
\end{theorem}

Throughout this paper, we always choose the constant $z_k$ of
Burkholder's inequality to denote the optimal constant
in Burkholder's $L^k(\P)$-inequality for continuous square-integrable martingales.
The precise value of $z_k$ involves the zeros of Hermite polynomials;
see  Davis \cite{Davis}. 

By the It\^o isometry, $z_2=1$. Carlen and Kree \cite[Appendix]{CK} have shown 
that $z_k \le 2\sqrt{k}$ for all $k \ge 2$, and moreover $z_k=(2+o(1))\sqrt{k}$
as $k\to\infty$.

%Throughout this paper, $z_{\nu}$ denotes the optimal constant
%in Burkholder's $L^\nu$-inequality for continuous square-integrable martingales; 
%its precise value has been computed by B. Davis \cite{Davis}. 
%By the It\^o isometry, $z_2=1$. Carlen and Kree \cite[Appendix]{CK} have shown 
%that $z_{\nu}\le 2\sqrt{\nu}$ for all $\nu\ge 2$, and moreover $z_{\nu}=(2+o(1))\sqrt{\nu}$
%as $\nu\to\infty$.

\section{Proof of Theorem \ref{th:main}: upper bound}

In this section we prove that $\overline\lambda(\nu)<\infty$
for all $\nu\in[2\,,\infty)$. 

If $\nu_1\le\nu_2$ are both in $[1\,,\infty)$, then by Jensen's inequality, 
\begin{equation}
	\limsup_{t\to\infty} \frac 1t \sup_{|x|\ge\alpha t}\ln
	\E\left(|u_t(x)|^{\nu_1}\right) \leq \frac{\nu_1}{\nu_2} 
	\limsup_{t\to\infty} \frac 1t \sup_{|x|\ge\alpha t}\ln
	\E\left(|u_t(x)|^{\nu_2}\right).
\end{equation}
This leads to the inclusion
\begin{equation}\begin{split}
	&\left\{\alpha>0:\
		\limsup_{t\to\infty} \frac 1t \sup_{|x|\ge\alpha t}\ln
		\E\left(|u_t(x)|^{\nu_2}\right)<0\right\} \\
	&\hskip1.5in\subseteq \left\{\alpha>0:\
		\limsup_{t\to\infty} \frac 1t \sup_{|x|\ge\alpha t}\ln
		\E\left(|u_t(x)|^{\nu_1}\right)<0\right\},
\end{split}\end{equation}
and hence the inequality
$\overline\lambda(\nu_1) \leq \overline\lambda(\nu_2)$. 
Therefore, it suffices to prove the result in the case that $\nu$ is an even integer $\ge 2$.
Our method is motivated strongly by ideas of Lunardi \cite{Lunardi} on
optimal regularity of analytic semigroups.

Dalang's condition \eqref{Dalang} implies
that the L\'evy process $X$ has transition functions
$p_t(x)$ \cite[Lemma 8.1]{FKN}; that is, for all
measurable $f:\R\to\R_+$,
\begin{equation}
	(P_tf)(x) := \E f(X_t)= \int_{-\infty}^\infty p_t(z)f(z)\,\d z
	\qquad\text{for all $t>0$}.
\end{equation}
And Dalang's theory implies that the solution can be written in 
mild form, in the sense of Walsh \cite{Walsh}, as
\begin{equation}\label{mild}
	u_t(x) = (P_tu_0)(x) + \int_{[0,t]\times\R} p_{t-s}(y-x)
	\sigma(u_s(y))\,W(\d s\,\d y),
\end{equation}
where $\{P_t\}_{t\ge 0}$ denotes the semigroup associated to the
process $X$. Henceforth, we will be concerned solely with
the mild formulation of the solution, as given to us by \eqref{mild}.

%Let ${\rm Lip}_\sigma$ denote the optimal Lipschitz constant of $\sigma$.

The following implies part 1 of Theorem \ref{th:main} immediately.

\begin{proposition}\label{pr:UB}
	If $\sup_{x\in\R}|\e^{cx/2}u_0(x)|$ and
	$\E\exp(cX_1)$ are both finite for some $c\in\R$, then
	for every even integer $\nu\ge 2$
	and for all
	\begin{equation}\label{cond:UB}
		\beta> \ln\E\e^{cX_1}+\frac 12 \Upsilon^{-1}\left(
		\left(2 z_{\nu} {\rm Lip}_\sigma\right)^{-2} \right),
	\end{equation}
	there exists a finite constant $A_{\beta,\nu}$
	such that $\E ( |u_t(x)|^{\nu} ) \le A_{\beta,\nu}\exp(\beta t-cx)$,
	uniformly for all $t\ge 0$ and $x\in\R$.
\end{proposition}

Proposition \ref{pr:UB} will be proved in Section \ref{subsec:proof:UB}.

\begin{remark}
	The proof shows that we require only that $\sigma(0)=0$;
	the positivity of ${\rm L}_\sigma$---see \eqref{cond:sigma}---is
	not required for this portion.
	
\end{remark}

\begin{remark} \label{rem:estimates}
	Proposition \ref{pr:UB} can frequently be used to give an explicit bound
	on $\overline\lambda(\nu)$. For example, if $\E\e^{c|X_1|}<\infty$
	for all $c\in\R$ and $u_0$ has compact support, then
	Proposition \ref{pr:UB} implies that
	\begin{equation}\label{LHS1}
		\limsup_{t\to\infty}\frac1t \ln\sup_{|x|\ge\alpha t}\E\left(|u_t(x)|^{\nu}\right)
		\le -\Lambda(\alpha)+\frac 12 \Upsilon^{-1}\left(
		\left(2 z_{\nu} {\rm Lip}_\sigma\right)^{-2}\right),
	\end{equation}
	where $\Lambda(\alpha):=\sup_{c\in\R}(\alpha c-\ln\E\e^{cX_1})$
	is the Legendre transformation of the logarithmic moment-generating function
	of $X_1$; see, for example,
	Dembo and Zeitouni \cite{DZ}. 
	Thus, the left-hand side of \eqref{LHS1} is negative as soon as
	$\Lambda(\alpha)>\frac 12 \Upsilon^{-1} ((2 z_{\nu} {\rm Lip}_\sigma)^{-2} )$, and hence
	\begin{equation}
		\overline{\lambda}(\nu) \le \inf\left\{ \alpha>0:\
		\Lambda(\alpha) > \frac 12 \Upsilon^{-1}\left(
		\left(2 z_{\nu} {\rm Lip}_\sigma\right)^{-2}\right)\right\}.
	\end{equation}
	We do not know how to obtain useful explicit lower bounds for 
	$\underline\lambda(\nu)$ in general. However,
	when $\mathcal{L} f= \frac{\kappa}{2} f''$, Theorem \ref{th:heat_exact}
	contains more precise bounds for both indices
	$\underline{\lambda}(2)$ and $\overline{\lambda}(2)$.
	
\end{remark}

\subsection{Stochastic weighted Young inequalities} \label{sec:young}

Proposition \ref{pr:UB} is based on general principles that 
might be of independent interest. These results will also be used 
in Section \ref{sec:wave} to study a family of
hyperbolic SPDEs. Throughout this subsection, $\Gamma_t(x)$ 
defines a nonrandom measurable function on $(0\,,\infty)\times\R$,
and $Z$ a predictable random field \cite[p.\ 292]{Walsh}. 

Consider the stochastic convolution
\begin{equation} \label{eq:conv}
	(\Gamma*Z\dot W)_t(x):=\int_{[0,t]\times\R} \Gamma_{t-s}(y-x)
	Z_s(y)\, W(\d s\,\d y),
\end{equation}
provided that it is defined in the sense of Walsh \cite[Theorem 2.5]{Walsh}. According
to the theory of Walsh, when
it is defined,
$\Gamma*Z\dot W$ defines a predictable random field. We 
study its $L^{\nu}(\P)$ norm next.

\begin{lemma}\label{lem:Young:conv:1}
	For all even integers $\nu\ge 2$, $t\ge 0$, and $x\in\R$,
	\begin{equation}
		\left\| (\Gamma*Z\dot W)_t(x)\right\|_{\nu}
		\le z_{\nu} \left(\int_{[0,t]\times\R}
		\Gamma^2_{t-s}(y-x)\left\|Z_s(y)\right\|_{\nu}^2\,\d s\,\d y
		\right)^{1/2},
	\end{equation}
	where $z_{\nu}$ was defined in Theorem \ref{th:BDG}.
\end{lemma}
\begin{proof}
	For fixed $t>0$ and $x\in\R$, we apply Burkholder's inequality 
	(Theorem \ref{th:BDG}) to the martingale 
	\begin{equation}
	r \mapsto \int_{[0,r]\times\R} \Gamma_{t-s}(y-x)
	Z_s(y)\, W(\d s\,\d y),
	\end{equation} 
	which has quadratic variation given by 
	\begin{equation}
	r \mapsto \int_{[0,r]\times\R}
			\Gamma^2_{t-s}(y-x) Z_s(y)^2\,\d s\,\d y.
	\end{equation} 
	We let $r=t$ to obtain
	\begin{equation}\begin{split}
		&\left\| (\Gamma*Z\dot W)_t(x) \right\|_{\nu}^{\nu}
			\le z_{\nu}^{\nu} \E\left(\left|\int_{[0,t]\times\R}
			\Gamma^2_{t-s}(y-x) Z_s(y)^2\,\d s\,\d y \right|^{\nu/2}\right)\\
		&\hskip.8in=z_{\nu}^{\nu} \E\left(\int_{([0,t]\times\R)^{\nu/2}}
			\prod_{j=1}^{\nu/2}
			\Gamma^2_{t-s_j}(y_j-x)|Z_{s_j}(y_j)|^2\,\d\bm{s}\,\d\bm{y}
			\right).
	\end{split}\end{equation}
	The generalized H\"older inequality implies that
	\begin{equation}
		\E\left(\prod_{j=1}^{\nu/2} |Z_{s_j}(y_j)|^2\right)
		\le \prod_{j=1}^{\nu/2}\left\|Z_{s_j}(y_j)\right\|_{\nu}^2,
	\end{equation}
	and the result follows.
 \end{proof}

We say that $\vartheta:\R\to\R_+$ is
a \emph{weight} when $\vartheta$ is measurable
and 
\begin{equation}\label{theta:submult}
	\vartheta(a+b)\le\vartheta(a)\vartheta(b)
	\qquad\text{for all $a,b\in\R$}.
\end{equation}
As usual, the weighted $L^2$-space $L^2_\vartheta(\R)$ 
denotes the collection of all measurable
functions $h:\R\to\R$ such that $\|h\|_{L^2_\vartheta(\R)}<\infty$, where
\begin{equation}
	\|h\|_{L^2_\vartheta(\R)}^2:=
	\int_{-\infty}^\infty |h(x)|^2\,\vartheta(x)\,\d x.
\end{equation}

Define, for all predictable processes $v$,
$\nu\in[1\,,\infty)$, and $\beta>0$,
\begin{equation}
	\mathcal{N}_{\beta,\nu,\vartheta} (v) := \left[
	\sup_{t\ge 0}\sup_{x\in\R}
	\e^{-\beta t} \vartheta(x)\left\| v_t(x)\right\|_{\nu}^2 \right]^{1/2}.
\end{equation}

\begin{proposition}[A stochastic Young inequality]\label{pr:young:1}
	For all weights $\vartheta$, all $\beta>0$, and all even integers $\nu\ge 2$,
	\begin{equation}
		\mathcal{N}_{\beta,\nu,\vartheta}(\Gamma*Z\dot W) \le
		z_{\nu}  \left(
		\int_0^\infty 
		\e^{-\beta t}\left\| \Gamma_t\right\|_{L^2_\vartheta(\R)}^2\,\d t\right)^{1/2}
		\cdot\mathcal{N}_{\beta,\nu,\vartheta}(Z).
	\end{equation}
\end{proposition}

\begin{proof}
	We apply Lemma \ref{lem:Young:conv:1} together with
	\eqref{theta:submult} to find that
	\begin{equation}\begin{split}
		&\e^{-\beta t}\vartheta(x)\left\| (\Gamma*Z\dot W)_t(x)\right\|_{\nu}^2\\
		&\le z_{\nu}^2\int_{[0,t]\times\R}\e^{-\beta (t-s)}
			\vartheta(y-x)\Gamma^2_{t-s}(y-x)\,\e^{-\beta s} \vartheta(y)
			\left\| Z_s(y)\right\|_{\nu}^2\,\d s\,\d y\\
		&\le z_{\nu}^2 \left|
			\mathcal{N}_{\beta,\nu,\vartheta}(Z)\right|^2\cdot\int_{[0,t]\times\R}
			\e^{-\beta r}\vartheta(z)\Gamma^2_r(z)\,\d r\,\d z.
	\end{split}\end{equation}
	The proposition follows from optimizing this expression over all $t\ge 0$
	and $x\in\R$.
 \end{proof}

\begin{proposition}\label{pr:conv:UB}
	If $\E\exp(cX_1)<\infty$ for some $c\in\R$, then for all predictable random fields
	$Z$, all $\beta>\ln\E\e^{cX_1}$, and all even integers $\nu\ge 2$,
	\begin{equation}
		\mathcal{N}_{\beta,\nu,\vartheta_c}(p*Z\dot W) \le
		z_{\nu}\left(2 \Upsilon\left(
		2 \beta- 2\ln\E\e^{cX_1}\right)\right)^{1/2}\cdot\mathcal{N}_{\beta,\nu,
		\vartheta_c}(Z),
	\end{equation}
	where $\vartheta_c(x):=\exp(cx)$.
\end{proposition}

\begin{proof}
	If $\vartheta$ is an arbitrary weight, then $\|p_t\|_{L^2_\vartheta(\R)}^2
	\le \sup_{z\in\R}p_t(z)\cdot \E\,\vartheta(X_t)$.
	According to the inversion formula,
	\begin{equation}\label{inversion}
		\sup_{z\in\R}p_t(z) \le \frac{1}{2\pi}\int_{-\infty}^\infty
		\e^{-t\Re\Psi(\xi)}\,\d\xi,
	\end{equation}
	whence
	\begin{equation}\begin{split}
		\int_0^\infty \e^{-\beta t}\|p_t\|_{L^2_\vartheta(\R)}^2\,\d t
			&\le \frac{1}{2\pi}\int_{-\infty}^\infty\d\xi
			\int_0^\infty\d t \ \e^{-t(\beta+\Re\Psi(\xi))}\E\,\vartheta(X_t).
	\end{split}\end{equation}
	The preceding is valid for all weights
	$\vartheta$. Now consider the following special case of $\vartheta:=\vartheta_c$.
	Clearly, this is a weight and,  in addition, by standard facts about
	L\'evy processes,
	\begin{equation}\label{E:theta:c}
		\E\,\vartheta_c(X_t) =\left( \E\e^{cX_1}\right)^t.
	\end{equation}
	Consequently, for all $\beta>M(c):= \ln\E\e^{cX_1}$,
	\begin{equation}\begin{split}
		\int_0^\infty \e^{-\beta t}\|p_t\|_{L^2_{\vartheta_c}(\R)}^2\,\d t
			&\le \frac{1}{2\pi}\int_{-\infty}^\infty\d\xi
			\int_0^\infty\d t \ \e^{-t(\beta+\Re\Psi(\xi)-M(c))}\\
		&=2\Upsilon(2\beta-2M(c)).
	\end{split}\end{equation}
	Proposition \ref{pr:young:1} completes the proof.
 \end{proof}

\begin{lemma}\label{u_0:young}
	For all  weights $\vartheta$, all
	$\beta>0$, and all even integers $\nu\ge 2$,
	\begin{equation}\label{eq:u_0:young}
		\mathcal{N}_{\beta,\nu,\vartheta}(P_\bullet u_0) \le
		 \mathcal{N}_{\beta,\nu,\vartheta}(u_0)\cdot \sup_{t\ge 0}
		\left(\e^{-\beta t}\E\,\vartheta(X_t)\right)^{1/2},
	\end{equation}
	where $P_\bullet u_0$ stands for the function $t \mapsto (P_t u_0)(x)$.	In particular, if $\E\e^{cX_1}<\infty$ for some $c\in\R$, then for all
	$\beta>\ln\E\e^{cX_1}$,
	\begin{equation}
		\mathcal{N}_{\beta,\nu,\vartheta_c}(P_\bullet u_0)
		\le \mathcal{N}_{\beta,\nu,\vartheta_c}(u_0).
	\end{equation}
\end{lemma}

%\begin{remark}
%	Since $u_0$ and $P_tu_0$ are nonrandom, it is clear that
%	\begin{equation}\begin{split}
%		\mathcal{N}_{\beta,\nu,\vartheta}(P_\bullet u_0)
%			&=\sup_{\substack{t\ge 0\\
%			x\in\R}}\left[|\vartheta(x)|^{1/2}|(P_tu_0)(x)|\right],\\
%		\mathcal{N}_{\beta,\nu,\vartheta}( u_0) &=
%			\sup_{x\in\R}|\vartheta(x)|^{1/2}|u_0(x)|.
%	\end{split}\end{equation}
%	But we have found it more helpful to write 
%	the norms abstractly.
%	
%\end{remark}

\begin{proof}
	%Without loss of generality $u_0(x)\ge 0$ for all $x$; otherwise,
	%we consider $|u_0(x)|$ in place of $u_0(x)$ in the following.
	
	Thanks to \eqref{theta:submult},
	\begin{align}\nonumber
		|\vartheta(x)|^{1/2} (P_tu_0)(x) &\le 
			\int_{-\infty}^\infty |\vartheta(y-x)|^{1/2}p_t(y-x)
			|\vartheta(y)|^{1/2}u_0(y)\,\d y\\
		&\le  \sup_{y\in\R}\left[ |\vartheta(y)|^{1/2} u_0(y)\right]\cdot
			\E\left(|\vartheta(X_t)|^{1/2}\right).
	\end{align}
	This and the Cauchy--Schwarz inequality together
	imply \eqref{eq:u_0:young}, and
	the remainder of the lemma follows from \eqref{E:theta:c}.
 \end{proof}

\subsection{Proof of Proposition \ref{pr:UB}} \label{subsec:proof:UB}
We begin by studying the Picard-scheme approximation to the solution $u$.
Namely, let $u^{(0)}_t(x):=u_0(x)$, and then define iteratively
\begin{equation}
	u_t^{(n+1)}(x) := (P_tu_0)(x) + 
	\left(p*\left(\sigma\circ u^{(n)}\right)\dot W\right)_t(x),
\end{equation}
for $t>0$, $x\in\R$, and $n\ge 0$, where the 
stochastic convolution is defined in \eqref{eq:conv}. Clearly,
\begin{equation} \label{eqn:picard}
	\left\| u^{(n+1)}_t(x)\right\|_{\nu} \le |(P_tu_0)(x)| +	
	\left\| \left(p*\left(\sigma\circ u^{(n)}\right)\dot W\right)_t(x)\right\|_{\nu},
\end{equation}
whence for all $\beta>\ln\E\e^{cX_1}$,
\begin{equation} \label{eqn:picard_2}
	\mathcal{N}_{\beta,\nu,\vartheta_c}\left( u^{(n+1)}\right)
	\le \mathcal{N}_{\beta,\nu,\vartheta_c}(u_0) + 
	z_{\nu}{\rm Lip}_\sigma
	\mathcal{T}^{1/2}
	\cdot\mathcal{N}_{\beta,\nu,\vartheta_c}
	\left(u^{(n)}\right),
\end{equation}
where $\mathcal{T}:= 2\Upsilon (
2\beta-2\ln\E\e^{cX_1})$;
see Proposition \ref{pr:conv:UB} and Lemma \ref{u_0:young}.
Condition \eqref{cond:UB} is equivalent to the  inequality
$z_{\nu}^2{\rm Lip}_\sigma^2 \mathcal{T}<1$.
Therefore, it follows from iteration that the quantity
$\mathcal{N}_{\beta,\nu,\vartheta_c}(u^{(n+1)})$ is bounded uniformly in $n$,
for this choice of $\beta$. Dalang's theory  \cite[Theorem 13 and its proof]{Dalang}
tells us that $\lim_{n\to\infty} u^{(n)}_t(x)=u_t(x)$
in probability for all $t\ge 0$ and $x\in\R$. Therefore,
Fatou's lemma
implies that $\mathcal{N}_{\beta,\nu,\vartheta_c}(u)<\infty$
when $\beta>\ln\E\e^{cX_1}$. This completes the
proof of Proposition \ref{pr:UB} [and hence part 1 of Theorem \ref{th:main}].

\section{Proof of Theorem \ref{th:main}: lower bound}\label{sec:2nd}

Our present, and final, goal is to prove that for all $\nu\in[2\,,\infty)$,  whenever 
$0<\alpha$ is sufficiently small,
$\limsup_{t\to\infty} t^{-1} \sup_{|x|>\alpha t}\ln 
\| u_t(x)\|_{\nu} >0$.
By Jensen's inequality, it suffices to prove this in the case that $\nu=2$.
We will borrow liberally several localization
ideas from two related papers by Mueller \cite{Mueller}
and Mueller and Perkins \cite{MuellerPerkins}.

Define, for all predictable random fields $v$, and $\alpha,\beta>0$,
\begin{equation} \label{def:norm_M}
	\mathcal{M}_{\alpha,\beta} (v) := \left[
	\int_0^\infty \e^{-\beta t}\,\d t\int_{\substack{x\in\R:\\
	|x|\ge\alpha t}} 
	\,\d x\ \left\| v_t(x)\right\|_2^2 \right]^{1/2}.
\end{equation}
Thus, $\{\mathcal{M}_{\alpha,\beta}\}_{\alpha,\beta>0}$ defines a family of
norms on the family of predictable random fields.

\begin{proposition}\label{pr:LB:apriori}
	If  $\E\vert X_1\vert<\infty$,
	then $\mathcal{M}_{\alpha,\beta}(u)=\infty$ for
	all sufficiently small $\alpha,\beta>0$.
\end{proposition}

\begin{proof}
	Thanks to \eqref{mild} and the It\^o isometry for stochastic integrals, 
	\begin{equation}\label{L2LB}
		\left\|u_t(x)\right\|_2^2 \ge \left|(P_tu_0)(x)\right|^2
		+{\rm L}_\sigma^2 \cdot\int_0^t\d s\int_{-\infty}^\infty\d y\
		\left|p_{t-s}(y-x)\right|^2\left\|u_s(y)\right\|_2^2.
	\end{equation}
	Let us define
	\begin{equation}
		\mathcal{M}^{+}_{\alpha,\beta} (v) := \left[
		\int_0^\infty \e^{-\beta t}\,\d t\int_{\substack{x\in\R:\\
		x\ge\alpha t}} 
		\,\d x\ \left\| v_t(x)\right\|_2^2 \right]^{1/2},
	\end{equation}
	and
	\begin{equation}
		\mathcal{M}^{-}_{\alpha,\beta} (v) := \left[
		\int_0^\infty \e^{-\beta t}\,\d t\int_{\substack{x\in\R:\\
		x\le-\alpha t}} 
		\,\d x\ \left\| v_t(x)\right\|_2^2 \right]^{1/2}.
	\end{equation}
	
	If $x,y\in\R$ and $0\le s\le t$, then the triangle inequality
	implies that
	\begin{equation}\label{eq:Ind:Ineq}
		\1_{[\alpha t,\infty)}(x) \ge \1_{[\alpha(t-s),\infty)}(x-y)
		\cdot\1_{[\alpha s,\infty)}(y).
	\end{equation}
	For all $r\ge 0$, let
	\begin{equation}
		T^{+}_\alpha(r):=\int_{\substack{z\in\R:\\z\ge\alpha r}}
		\left| p_r(z)\right|^2\,\d z,\quad
		T^{-}_\alpha(r):=\int_{\substack{z\in\R:\\z\le-\alpha r}}
		\left| p_r(z)\right|^2\,\d z,
	\end{equation}
	and
		\begin{equation}
		S^{+}_\alpha(r) :=\int_{\substack{y\in\R:\\y\ge\alpha r}}
		\|u_r(y)\|_2^2\,\d y,\quad
		S^{-}_\alpha(r) :=\int_{\substack{y\in\R:\\y\le-\alpha r}}
		\|u_r(y)\|_2^2\,\d y.
	\end{equation}
	According to \eqref{eq:Ind:Ineq},
	\begin{equation}\label{conv}
		\int_{x\ge\alpha t}\|u_t(x)\|_2^2\,\d x
		\ge \int_{x\ge\alpha t}\left| (P_t u_0)(x)\right|^2\,\d x
		+{\rm L}_\sigma^2\cdot(T^{-}_\alpha*S^{+}_\alpha)(t),
	\end{equation}
	where ``$*$'' denotes the usual convolution on $\R_+$.
	
	We multiply both sides of
	\eqref{conv} by $\exp(-\beta t)$ and integrate
	$[\d t]$ to find
	\begin{equation}\begin{split} \label{Mplus}
		\left|\mathcal{M}^{+}_{\alpha,\beta}(u) \right|^2&\ge 
			\left|\mathcal{M}^{+}_{\alpha,\beta}(P_\bullet u_0)\right|^2+
			{\rm L}_\sigma^2 \cdot \widetilde{T}^{-}_\alpha(\beta)
			\widetilde{S}^{+}_\alpha(\beta)\\
		&=\left|\mathcal{M}^{+}_{\alpha,\beta}(P_\bullet u_0)\right|^2+
			{\rm L}_\sigma^2 \cdot \widetilde{T}^{-}_\alpha(\beta)
			\left| \mathcal{M}^{+}_{\alpha,\beta}(u) \right|^2,
	\end{split}\end{equation}
	where $\widetilde H(\beta):=\int_0^\infty\exp(-\beta t)H(t)\,\d t$
	defines the Laplace transform of $H$ for every measurable function
	$H:\R_+\to\R_+$. Also, we can apply a similar argument, run on the 
	negative half of the real line, to deduce that
	\begin{equation} \label{Mminus}
		\left|\mathcal{M}^{-}_{\alpha,\beta}(u) \right|^2 \ge 
		 \left|\mathcal{M}^{-}_{\alpha,\beta}(P_\bullet u_0)\right|^2+
			{\rm L}_\sigma^2 \cdot \widetilde{T}^{+}_\alpha(\beta)
			\left| \mathcal{M}^{-}_{\alpha,\beta}(u) \right|^2.
	\end{equation}
	Next we add the inequalities \eqref{Mplus} and \eqref{Mminus}:
	Because $\{X_t\}_{t\ge 0}$ is symmetric,
	$\tilde{T}^{+}_{\alpha}(\beta) = \tilde{T}^{-}_{\alpha}(\beta)$;
	and it is easy to see that
	$\mathcal{M}_{\alpha,\beta}(u)^2 = 
	\mathcal{M}^{+}_{\alpha,\beta}(u)^2 + \mathcal{M}^{-}_{\alpha,\beta}(u)^2$.
	Therefore, we can conclude that
	\begin{equation} \label{M}
		\left|\mathcal{M}_{\alpha,\beta}(u) \right|^2 \ge 
		 \left|\mathcal{M}_{\alpha,\beta}(P_\bullet u_0)\right|^2+
		{\rm L}_\sigma^2 \cdot \widetilde{T}^{+}_\alpha(\beta)
		\left| \mathcal{M}_{\alpha,\beta}(u) \right|^2.
	\end{equation}
	
	Next we may observe that
	\begin{equation}
		|\mathcal{M}_{\alpha,\beta}(P_\bullet u_0)|>0.
	\end{equation}
	This holds because $u_0\ge 0$, $u_0>0$ on a set of positive measure,
	and $u_0$ is lower semicontinuous. Indeed, if it were not
	so, then $\int_{|x|\ge\alpha t} (P_tu_0)(x)\,\d x=0$ for almost all,
	hence all, $t>0$. But then we would let $t\to 0$ to deduce from this and
	Fatou's lemma that
	$\int_{-\infty}^\infty u_0(x)\,\d x=0$, which is a contradiction.
	
	The preceding development implies the following:
	\begin{equation}\label{ifthen}
		\text{If	} \mathcal{M}_{\alpha,\beta}(u)<\infty,
		\text{ then } \widetilde{T}^{+}_\alpha(\beta) < {\rm L}_\sigma^{-2}.
	\end{equation}
	
	The symmetric L\'evy process $X$ is recurrent iff
	\begin{equation}\label{inf}
		\Upsilon(0^+)=\infty.
	\end{equation}
	See, for example, Port and Stone \cite[\S16]{PortStone}.
	Therefore it remains to prove
	that the conditions of Theorem \ref{th:main}
	imply \eqref{inf}.
	
	The discrete-time process $\{X_n\}_{n=1}^\infty$
	is a one-dimensional mean-zero [in fact symmetric] random walk, which is necessarily
	recurrent thanks to the Chung--Fuchs theorem \cite{ChungFuchs}. 
	Consequently, the L\'evy process $\{X_t\}_{t\ge 0}$ is recurrent as well.
	Thanks to the preceding paragraph, \eqref{inf} holds.
	
	By the monotone convergence theorem, 
	\begin{equation}
		\lim_{\alpha\downarrow 0}
		\widetilde{T}_\alpha(\beta)= \frac 12 \int_0^\infty 
		\e^{-\beta t}\|p_t\|_{L^2(\R)}^2\,\d t
		=\frac 12 \Upsilon(\beta)\qquad\text{for all $\beta>0$}.
	\end{equation}
	[The second identity follows from Plancherel's theorem.] Let $\beta\downarrow 0$
	and appeal to \eqref{inf}
	to conclude that $\widetilde{T}^+_\alpha(\beta)>{\rm L}_\sigma^{-2}$
	for all sufficiently-small positive $\alpha$ and $\beta$. In light
	of \eqref{ifthen}, this completes our demonstration.
 \end{proof}

\begin{proof}[Proof of Part 2 of Theorem \ref{th:main}]
	
	Choose and fix $\alpha$ and $\beta$ positive, but so small that
	$\mathcal{M}_{\alpha,\beta}(u)=\infty$ [Proposition \ref{pr:LB:apriori}].
	According to Proposition \ref{pr:UB},  for all fixed $\alpha'>0$,
	\[
		\int_0^\infty\e^{-\beta t}\,\d t\int_{|x|\ge\alpha' t}\d x\
		\|u_t(x)\|_2^2 
		\le A_{\alpha',2}\int_0^\infty 
		\e^{(\beta'-\beta)t}\,\d t\int_{|x|\ge\alpha' t} \d x\ \e^{-c|x|},
	\]
	provided that $\beta'$ [in place of the variable $\beta$ there]
	satisfies \eqref{cond:UB} with $\pm c$ [in place of the variable $c$ there].
	We choose and fix $\beta'$ so large that the condition \eqref{cond:UB}
	is satisfied for $\beta'$. Then,  choose and fix $\alpha'$ so large
	that the right-most integral in the preceding display is finite. 
	Since $\mathcal{M}_{\alpha,\beta}(u)=\infty$,
	it follows from the preceding that
	\begin{equation}
		\int_0^\infty \e^{-\beta t}\,\d t\int_{\alpha t\le |x|\le \alpha' t}
		\d x\ \E \left(|u_t(x)|^2 \right)=\infty.
	\end{equation}
	Consequently,
	\begin{equation}
		\int_0^\infty t \e^{-\beta t} \sup_{|x|\ge\alpha t}\E \left(
		|u_t(x)|^2 \right)\,\d t=\infty,
	\end{equation}
	whence
	\begin{equation}
		\limsup_{t\to\infty}t^{-1} \sup_{|x| \ge \alpha t} \ln
		\E \left(|u_t(x)|^2 \right)\ge\beta>0
	\end{equation}
	for the present choice of $\alpha$ and $\beta$. This implies that
	$\underline{\lambda}(2)\ge\alpha>0$.
 \end{proof}

\begin{remark}
	Theorem \ref{th:main} requires less than 
	the symmetry of the L\'evy process $\{X_t\}_{t\ge 0}$.
	For instance, our proof continues to work
	provided that there exist finite and positive constants $c_1$ and $c_2$ such that
	\begin{equation}
		c_1 T^{-}_{\alpha}(r)\leq T^{+}_{\alpha}(r) \leq c_2T^{-}_{\alpha}(r),
	\end{equation}
	simultaneously for all $\alpha>0$ and $r \geq 0$.

\end{remark}

\section{Proof of Theorem \ref{th:heat_exact}} \label{sec:heat_exact}

Throughout the proof, we choose and fix some $\kappa>0$.
Thus, the operator $\mathcal{L}f = \frac{\kappa}{2} f''$ 
is the generator of a L\'evy process given by 
$X_t = \sqrt{\kappa} \, B_t$, where $\{B_t\}_{t\ge 0}$ 
is a Brownian motion, and Theorem \ref{th:main} obviously applies in this case.
We now would like to prove the second claim of Theorem \ref{th:heat_exact}. 
We proceed as we did for Theorem \ref{th:main}, and 
divide the proof in two parts: One part is concerned with an upper bound 
for $\overline{\lambda}(2)$; and the other deals with a lower bound on $\underline{\lambda}(2)$.

\subsection{Upper bound}

In order to obtain an upper estimate for $\overline{\lambda}(2)$, 
we could follow the procedure outlined in Remark \ref{rem:estimates}. 
But this turns out to be not optimal. 
In the case of Theorem \ref{th:heat_exact}, we know explicitly the 
transition functions $p^{(\kappa)}_t$:
\begin{equation} \label{density}
	p^{(\kappa)}_t(x) = \frac{1}{\sqrt{2 \pi \kappa t}} \exp\left(-\frac{x^2}{2 \kappa t}\right).
\end{equation}
Therefore, we can use \eqref{density} directly and make exact computations
in order to improve on the general bounds of Remark \ref{rem:estimates}. 
We first prove the following; it sharpens Proposition \ref{pr:UB} in the present setting.

\begin{proposition} \label{pr:UB_heat}
	If $\mathcal{L} f= \frac{\kappa}{2} f''$ 
	and $\sup_{x\in\R} |\e^{cx/2}u_0(x)|$ is finite for some $c\in\R$, then for every
	\begin{equation}\label{cond:UB_heat}
		\beta > \frac{\kappa c^2}{4} + \frac{{\rm Lip}_\sigma^4}{4\kappa},
	\end{equation}
	there exists a finite constant $A_{\beta}$ such that 
	$\E ( |u_t(x)|^2 ) \le A_{\beta}\exp(\beta t-cx)$, uniformly for all $t\ge 0$ and $x\in\R$.
\end{proposition}

	\begin{proof}
	We follow the proof of Proposition \ref{pr:UB}, but use Proposition \ref{pr:young:1},
	instead of Proposition \ref{pr:conv:UB}, in order to handle \eqref{eqn:picard} better. 
	Then, \eqref{eqn:picard_2} is replaced by
	\begin{eqnarray}
		\lefteqn{\mathcal{N}_{\beta,2,\vartheta_c}\left( u^{(n+1)}\right)} \\
		& \le & \mathcal{N}_{\beta,2,\vartheta_c}(u_0) + {\rm Lip}_\sigma \left(\int_0^\infty 
			\e^{-\beta t} 
			\left\| p^{(\kappa)}_t\right\|_{L^2_{\vartheta_c}(\R)}^2 
			\, \d t \right)^{1/2} \mathcal{N}_{\beta,2,\vartheta_c}
			\left(u^{(n)}\right). \nonumber
	\end{eqnarray}
	Next we complete the proof, in the same way we did for Proposition \ref{pr:UB}, 
	and deduce that there exists a constant $A_{\beta}$ such that 
	$\E(|u_t(x)|^2) \le A_{\beta} \exp(\beta t-cx)$ uniformly for all $t \ge 0$ and
	$x \in \R$, provided that $\beta$ is chosen to be large enough to satisfy
	\begin{equation} \label{opt_cond}
		{\rm Lip}_\sigma^2 \cdot \int_0^\infty 
		\e^{-\beta t} \left\| p^{(\kappa)}_t\right\|_{L^2_{\vartheta_c}(\R)}^2 \, \d t  < 1.
	\end{equation}
	Now we compute:
	\begin{equation}\begin{split}
		\left\| p^{(\kappa)}_t\right\|_{L^2_{\vartheta_c}(\R)}^2 & = 
			\frac{1}{2 \pi \kappa t} \int_{-\infty}^\infty 
			\exp\left(-\frac{x^2}{\kappa t} + cx\right) \, \d x \\
		& = \frac{1}{2 \sqrt{\pi \kappa t}} \exp\left(\frac{\kappa c^2 t}{4}\right).
	\end{split}\end{equation}
	Since $\int_{0}^{\infty} t^{-1/2}e^{- \beta t} \, \d t 
	= \sqrt{ \pi/\beta}$, we have the following for
	all $\beta > \kappa c^2/4$:
	\begin{equation}
		{\rm Lip}_\sigma^2 \cdot \int_0^\infty 
		\e^{-\beta t} \left\| p^{(\kappa)}_t\right\|_{L^2_{\vartheta_c}(\R)}^2 \, \d t 
		= \frac12 {\rm Lip}_\sigma^2\left(\kappa \beta - \dfrac{\kappa^2 c^2}{4}\right)^{-1/2}.
	\end{equation}
	And hence, \eqref{opt_cond} follows from \eqref{cond:UB_heat}.
	This proves Proposition \ref{pr:UB_heat}.
 \end{proof}

\begin{proof}[Proof of the upper bound in Theorem \ref{th:heat_exact}]
	If $u_0$ has compact support, then
	 the assumption of Proposition \ref{pr:UB_heat} is satisfied for all $c \in \R$.
	Consequently, 
	\begin{equation}
	\limsup_{t\to\infty} \frac 1t \ln\sup_{|x|\ge\alpha t}
	\E (|u_t(x)|^2 ) \le \beta - c\alpha,
	\end{equation} 
	and hence
	\begin{equation}
		\overline{\lambda}(2) \le \inf\left\{ \alpha>0:\ \beta - c\alpha < 0\right\} = \frac{\beta}{c}.
	\end{equation}
	This and \eqref{cond:UB_heat} together imply that
	\begin{equation}
		\overline{\lambda}(2) \le \inf_{c\in\R}\left(
		\frac{\kappa c}{4} + \frac{{\rm Lip}_\sigma^4}{4\kappa c}\right)=\frac{{\rm Lip}_\sigma^2}{2}.
	\end{equation}
	This concludes the proof of the upper bound.
 \end{proof}

\subsection{Lower bound}

We first prove the following refinement of Proposition \ref{pr:LB:apriori}.

\begin{proposition} \label{pr:heat:LB}
	If $\mathcal{L}f = \frac{\kappa}{2} f''$ and $\alpha$ and $\beta$ satisfy
	\begin{equation} \label{cond_alpha}
		\left( \alpha - \frac{{\rm L}_\sigma^2}{4 \pi} \right)^2 <
		\frac{{\rm L}_\sigma^4}{16 \pi^2} - \kappa \beta,
	\end{equation}
	then $\mathcal{M}_{\alpha,\beta}(u) = \infty$.
\end{proposition}

\begin{proof}
	In the case that we consider here, the L\'evy process is a scaled Brownian motion. 
	Hence, Proposition \ref{pr:LB:apriori} applies, and in accord with \eqref{ifthen}, 
	it suffices to prove the following:
	\begin{equation} \label{heat:LB:cond}
		\mathfrak{I}:=
		\int_0^\infty \e^{-\beta t} \left( \int_{\substack{z\in\R:\\z\ge\alpha t}} 
		\left| p^{(\kappa)}_t(z)\right|^2\,\d z \right) \d t > {\rm L}_\sigma^{-2}.
	\end{equation}
	
	Let $\bar\Phi(z):=(2\pi)^{-1/2}\int_z^\infty \exp(-\tau^2/2)\,\d\tau$
	for every $z\in\R$,
	then apply \eqref{density} and compute directly to find that
	\begin{equation}\begin{split}
		\mathfrak{I}& = \frac{1}{2 \sqrt{\pi \kappa}} \int_0^\infty \frac{\e^{-\beta t}}{\sqrt{t}} 
			\bar\Phi\left(\sqrt{\frac{2\alpha^2 t}{\kappa}} \right) \d t\\
		& = \frac{\alpha}{4 \pi \kappa} \int_0^\infty 
			\frac{\e^{-\alpha^2 t/\kappa}}{\sqrt{t}} 
			\left(\int_{0}^{t} \frac{\e^{-\beta s}}{\sqrt{s}} \, \d s\right) \d t,
	\end{split}\end{equation}
	after we integrate by parts.
	Since $\e^{-\beta s} \ge \e^{-\beta t}$ for $s \le t$, 
	\begin{equation}
		\mathfrak{I} 
		\ge \frac{\alpha}{2 \pi \kappa} \int_0^\infty 
		\exp\left( - \left(\beta  + \frac{\alpha^2}{\kappa} \right) t\right) \d t 
		= \frac{\alpha}{2 \pi(\beta \kappa + \alpha^2)}.
	\end{equation}
	Hence, \eqref{cond_alpha} implies \eqref{heat:LB:cond}, and hence
	the proposition.
 \end{proof}

\begin{remark}
We notice that condition \eqref{cond_alpha} is sufficient but not necessary. Indeed, as $\mathfrak{I}$ is decreasing in $\alpha$, only the upper bound implied by \eqref{cond_alpha} is relevant. Typically, \eqref{heat:LB:cond} is satisfied for $\alpha = 0$.
\end{remark}

\begin{proof}[Proof of the lower bound in Theorem \ref{th:heat_exact}]
	The second part of the proof of Theorem \ref{th:main} shows 
	that $\underline{\lambda}(2) \ge \alpha$, provided that we choose
	$\alpha$ and $\beta$ such that $\mathcal{M}_{\alpha,\beta}(u) = \infty$. 
	In accordance with \eqref{cond_alpha}, and after maximizing over 
	$\beta \le {\rm L}_\sigma^4/(16 \pi^2 \kappa)$---that is, making
	$\beta$ as small as possible---we obtain
	$\underline{\lambda}(2) \ge \alpha \ge {\rm L}_\sigma^2/2 \pi$.
	This concludes the proof of Theorem \ref{th:heat_exact}.
 \end{proof}

\section{A nonlinear stochastic wave equation} \label{sec:wave}

In this section, we study the nonlinear stochastic wave equation
\begin{equation} \label{eqn:wave}
	\frac{\partial^2}{\partial t^2}u_t(x) = \kappa^2 \left(\frac{\partial^2}{\partial x^2} u_t\right) (x) + 
	\sigma(u_t(x)) \frac{\partial^2}{\partial t \partial x} W(t,x) 
	\quad \text{for} \quad t > 0,\, x \in \R,
\end{equation}
where: (i) $\sigma:\R\to\R$ is Lipschitz continuous with 
Lipschitz constant ${\rm Lip}_\sigma$; (ii) $W$ is two-parameter
Brownian sheet, indexed by $(t\,,x)\in\R_+\times\R$; 
(iii) the initial function $u_0:\R\to\R_+$ and the initial derivative 
$v_0:\R\to\R$ are both in $L^\infty(\R)$;
and (iv) $\kappa>0$.
In the present one-dimensional setting, 
the nonlinear equation \eqref{eqn:wave} has been studied by Carmona and Nualart
\cite{CarmonaNualart} and Walsh \cite{Walsh}. 
There are also results available in the more delicate setting where
$x\in\R^d$ for $d>1$; see Conus and Dalang \cite{ConusDalang}, 
Dalang \cite{Dalang}, Dalang and Frangos \cite{DalangFrangos}, and 
Dalang and Mueller \cite{DalangMueller}.

It is well known that the fundamental solution for the wave equation in spatial dimension 1 is 
\begin{equation} \label{eq:green_wave}
	\Gamma_t(x) := \frac{1}{2} \1_{[-\kappa t,\kappa t]}(x)\qquad\text{for
	$t>0$ and $x\in\R$}.
\end{equation}
According to the theory of 
Dalang \cite{Dalang}, the stochastic wave equation
\eqref{eqn:wave} has an a.s.-unique mild solution. 
In the case that $u_0$ and $v_0$ are both constant functions,
Dalang and Mueller \cite{DalangMueller2} have shown that the solution
to \eqref{eqn:wave} is intermittent. 

In this section we will use the stochastic weighted Young inequalities of
Section \ref{sec:young} in order to deduce the weak intermittence of 
the solution to \eqref{eqn:wave} for nonconstant functions $u_0$ and $v_0$. 
And more significantly, when $u_0$ and $v_0$ have compact support,
we describe the precise rate at which
the farthest peaks can move away from the origin.

Here and throughout,
we assume that \eqref{cond:sigma} holds, and define
$\overline{\lambda}(\nu)$ and $\underline{\lambda}(\nu)$ as in \eqref{def:lambda_up} and \eqref{def:lambda_down}.

\begin{theorem} \label{thm:wave}
	If $u_0, v_0 \in \mathcal{D}_{\text{\it exp}}$, $u_0 > 0$ 
	on a set of positive measure and $v_0 \geq 0$, then $0 < \underline{\lambda}(\nu) \le 
	\overline{\lambda}(\nu) < \infty$ for all $\nu \in [2\,,\infty)$. 
	If, in addition, $u_0$ and $v_0$ have compact support, then 
	$\underline{\lambda}(\nu) = \overline{\lambda}(\nu) = \kappa$
	for all  $\nu \in [2\,,\infty)$.
\end{theorem}

Theorem \ref{thm:wave} implies the weak intermittence of the solution to
\eqref{eqn:wave}. And more significantly, it tells that when the
initial data have compact support, we have a sharp phase 
transition [$\underline{\lambda}(\nu) = \overline{\lambda}(\nu) = \kappa$]:
The solution has exponentially-large peaks inside
$[-\kappa t+o(t)\,,\kappa t+o(t)]$, and is exponentially small everywhere
outside $[-\kappa t+o(t)\,,\kappa t+o(t)]$. In particular, the farthest high
peaks of the solution travel at sharp linear speed $\pm\kappa t+o(t)$.
This speed corresponds to the speed of the traveling waves if we consider the 
deterministic equivalent of \eqref{eqn:wave} [say, when $\sigma\equiv 0$]. 
We emphasize that, contrary to what happens in the stochastic heat equation 
(Theorem \ref{th:heat_exact}), the growth behavior of the solution
to the stochastic wave equation
\eqref{eqn:wave} depends on the size of 
the noise (that is, the magnitude of $\kappa$), but not on the
growth rate of the nonlinearity $\sigma$.

\subsection{Proof of Theorem \ref{thm:wave}: upper bound.}

The proof of Theorem \ref{thm:wave} follows closely those of
Theorems \ref{th:main} and \ref{th:heat_exact}. 

We first show that $\overline{\lambda}(\nu) < \infty$. 
The solution to \eqref{eqn:wave} can be written in mild form, as
\begin{equation}
	u_t(x) = U^{(0)}_t(x) + V^{(0)}_t(x) + \int_{[0,t]\times\R} 
	\Gamma_{t-s}(y-x) \sigma(u_s(y))\, W(\d s \, \d y),
\end{equation}
where $U^{(0)}_t(x) = \frac{1}{2}(u_0(x + \kappa t) + 
u_0(x - \kappa t))$ and $V^{(0)}_t(x) = \frac{1}{2\kappa} 
\int_{x-\kappa t}^{x+\kappa t} v_0(y) \, \d y$.

The following Proposition implies immediately 
that $\overline{\lambda}(\nu) < \infty$ for $\nu \ge 2$.

\begin{proposition} \label{pr:UB_wave}
	Let $\nu \ge 2$ be an even integer,
	and assume that both $\sup_{x\in\R} |\e^{cx/2}u_0(x)|$ 
	and $\sup_{x\in\R} |\e^{cx/2}v_0(x)|$ are
	finite for some $c\in\R$. Then for every
	\begin{equation}\label{cond:UB_wave}
		\beta > \sqrt{\kappa^2 c^2 + \frac{z_{\nu}^2{\rm Lip}_\sigma^2}{2}},
	\end{equation}
	there exists a finite constant $A_{\beta}$ such that 
	$\E ( |u_t(x)|^{\nu} ) \le A_{\beta}\exp(\beta t-cx)$, 
	uniformly for all $t\ge 0$ and $x\in\R$.
\end{proposition}

In order to prove Proposition \ref{pr:UB_wave}, we will 
need the following Lemma. Let $\vartheta_c$ and 
$\mathcal{N}_{\beta,\nu,\vartheta}$ be defined as they were in 
Section \ref{sec:young}.

\begin{lemma}\label{u_0_v_0:young}
	For all $c \in \R$, $\beta > \kappa |c|/2$, and even integers $\nu \ge 2$,
	\begin{equation}\label{eq:u_0:young_wave}
		\mathcal{N}_{\beta,\nu,\vartheta_c}(U^{(0)}) 
		\le \mathcal{N}_{\beta,\nu,\vartheta_c}(u_0)
		\quad\text{and}\quad
		\mathcal{N}_{\beta,\nu,\vartheta_c}(V^{(0)}) \le 
		\frac{1}{\kappa c} \mathcal{N}_{\beta,\nu,\vartheta_c}(v_0).
	\end{equation}
\end{lemma}

\begin{proof}
	The first inequality of \eqref{eq:u_0:young_wave} follows from the definition of
	$U^{(0)}$. As regards the second, we have
	\begin{equation}\begin{split}
		\e^{cx/2} V^{(0)}_t(x) & \le  \left(\sup_{y \in \R} \e^{cy/2} v_0(y)\right) 
			\frac{\e^{cx/2}}{2\kappa} \int_{x-\kappa t}^{x+\kappa t} \e^{-cy/2} \, \d y \\
		& \le  \frac{\e^{c\kappa t/2}}{\kappa c}\left(\sup_{y \in \R} \e^{cy/2} v_0(y)\right).
	\end{split}\end{equation}
	Because $\beta > \kappa |c|/2$, this proves the lemma.
 \end{proof}

\begin{proof}[Proof of Proposition \ref{pr:UB_wave}]
	As in the proof of Proposition \ref{pr:UB_heat}, we apply a Picard-iteration scheme
	to approximate the solution $u$. Then, Lemma \ref{u_0_v_0:young} 
	and Proposition \ref{pr:young:1} yield
	\begin{align}\label{eqn:picard_wave}
		\mathcal{N}_{\beta,\nu,\vartheta_c}\left( u^{(n+1)}\right) & 
			\le  \mathcal{N}_{\beta,\nu,\vartheta_c}(u_0) + \frac{1}{\kappa c}
			\mathcal{N}_{\beta,\nu,\vartheta_c}(v_0) \\ \nonumber
		& \quad + \; z_{\nu} {\rm Lip}_\sigma \left(\int_{0}^{\infty} 
			\e^{-\beta t} \|\Gamma_t\|^2_{L^2_{\vartheta_c}(\R)} \, \d t 
			\right)^{\frac{1}{2}} \cdot \mathcal{N}_{\beta,\nu,\vartheta_c}
			\left(u^{(n)}\right).
	\end{align}
	A direct computation, using only \eqref{eq:green_wave}, shows that
	\begin{equation}
		\int_{0}^{\infty} \e^{-\beta t} \|\Gamma_t\|^2_{L^2_{\vartheta_c}(\R)}
		\, \d t < 
		(z_{\nu}{\rm Lip}_\sigma)^{-2}.
	\end{equation}
	And the same arguments that were used in the proof of Proposition \ref{pr:young:1}
	can be used to deduce from this bound that
	$\mathcal{N}_{\beta,\nu,\vartheta_c}(u)$ is finite.
	Now we use \eqref{eq:green_wave} in order to see that this condition is equivalent to
	\eqref{cond:UB_wave}.
	This concludes the proof of Proposition \ref{pr:UB_wave}.
 \end{proof}

\begin{proof}[Proof of the upper bound in Theorem \ref{thm:wave}]
	Proposition \ref{pr:UB_wave} implies the upper bound $\overline{\lambda}(\nu) < \infty$. 
	Now suppose $u_0$ and $v_0$ have compact support.
	In that case, $c$ is an arbitrary real number. 
	And similar arguments as in the proof of the upper bound of 
	Theorem \ref{th:heat_exact} imply that $\overline{\lambda}(\nu) \le \beta/c$.
	Together with \eqref{cond:UB_heat}, this leads to the following estimate:
	\begin{equation}
		\overline{\lambda}(\nu) \le \inf_{c\in\R}
		\sqrt{\kappa^2 + \frac{2 z_{\nu}^2{\rm Lip}_\sigma^2}{c^2}}=\kappa.
	\end{equation}
	This proves half of the theorem.
 \end{proof}

\subsection{Proof of Theorem \ref{thm:wave}: lower bound}

The following proposition implies
the requisite bound for the second half of the proof of
Theorem \ref{thm:wave}; namely, that $\underline{\lambda}(\nu) > 0$ for $\nu \ge 2$. 
Let $\mathcal{M}_{\alpha,\beta}$ be defined as in \eqref{def:norm_M}.

\begin{proposition} \label{pr:wave:LB}
	$\mathcal{M}_{\alpha,\beta}(u) = \infty$ provided that 
	\begin{equation} \label{cond:alpha_wave}
		0 < \alpha < \kappa - \frac{4 \beta^2}{{\rm L}_\sigma^2}.
	\end{equation}
\end{proposition}

\begin{proof}
	Similar arguments as in the proof of Proposition \ref{pr:LB:apriori} show that
	\begin{equation}
		|\mathcal{M}_{\alpha,\beta}(u)|^2 \ge 
		|\mathcal{M}_{\alpha,\beta}(U^{(0)}+V^{(0)})|^2 + 
		{\rm L}_\sigma^2 \cdot \tilde{T}^+_{\alpha}(\beta) |\mathcal{M}_{\alpha,\beta}(u)|^2,
	\end{equation}
	where $\tilde{T}^+_{\alpha}(\beta)$ denotes the Laplace transform of 
	$T^+_\alpha(r) := \int_{z\ge\alpha r} | \Gamma_r(z)|^2\,\d z$. Since 
	$u_0 > 0$ on a set of positive measure and $v_0 \geq 0$, we have 
	$|\mathcal{M}_{\alpha,\beta}(U^{(0)} + V^{(0)})| > 0$. This shows that 
	if ${\rm L}_\sigma^2 \cdot \tilde{T}^+_{\alpha}(\beta) > 1$, then
	$\mathcal{M}_{\alpha,\beta}(u) = \infty$. A direct computation reveals that
	\begin{equation}
		\tilde{T}^+_{\alpha}(\beta) = \begin{cases}
			(\kappa - \alpha)/(4\beta^2) & \text{ if } \alpha \le \kappa,\\
			0 & \mbox{ otherwise. }
		\end{cases}
	\end{equation}
	Hence, $\mathcal{M}_{\alpha,\beta}(u) = \infty$ 
	if $\tilde{T}^+_{\alpha}(\beta) > {\rm Lip}_\sigma^{-2}$, and the latter condition is equivalent to
	\eqref{cond:alpha_wave}.
	Since we also want $\alpha > 0$, Proposition \ref{pr:wave:LB} follows.
 \end{proof}

\begin{proof}[Proof of the lower bound in Theorem \ref{thm:wave}]
	For every $\alpha$ such that $\mathcal{M}_{\alpha,\beta}(u) = \infty$, 
	we can apply the same arguments as in the proof of the lower bound of 
	Theorem \ref{th:main} in order to conclude 
	that $\underline{\lambda}(2) \ge \alpha > 0$. Now, Proposition \ref{pr:wave:LB} shows that
	$\underline{\lambda}(2) \ge \kappa - 4 \beta^2/{\rm L}_\sigma^2$
	for all $\beta>0$, whence
	$\underline{\lambda}(2) \ge \kappa$. 
	Jensen's inequality then shows that $\underline{\lambda}(\nu) \ge \kappa$ as well.
 \end{proof}

\begin{remark}
	The condition $v_0 \geq 0$ 
	is not necessary in Theorem \ref{thm:wave}. Indeed, the necessary 
	condition is
	\begin{equation}\label{eq:Mpos}
		\mathcal{M}_{\alpha,\beta}(U^{(0)} + V^{(0)}) > 0. 
	\end{equation}
	The easy-to-verify conditions on $u_0$ and $v_0$ [in Theorem \ref{thm:wave}] imply 
	\eqref{eq:Mpos}.
\end{remark}

\begin{small}{\bf Acknowledgements:}
	We would like to thank Le Chen, Robert Dalang, 
	and two anonymous referees for generously sharing with us
	their many corrections and comments which ultimately led
	to a much better corrected draft.
\end{small}

\begin{small}

\vskip.2in
\noindent\textbf{Daniel Conus} and \textbf{Davar Khoshnevisan}\\
\noindent Department of Mathematics, University of Utah,
		Salt Lake City, UT 84112-0090\\
\noindent\emph{Emails:} 		\texttt{conus@math.utah.edu},
	\texttt{davar@math.utah.edu}\\
\noindent\emph{URLs:}	\texttt{http://www.math.utah.edu/\~{}conus},
	\texttt{http://www.math.utah.edu/\~{}davar}\\
\end{small}

\end{document}